\documentclass[11pt]{article}
\usepackage{amsmath,amsfonts,calrsfs,amssymb,color}
   \large\normalsize
   
\begin{document}

    {\bf   Multiplicative controllability of  the  reaction-diffusion equation on a parallelepiped with finitely many zero hyperplanes.}
    

\bigskip
\par\noindent
A.Y. Khapalov  
\\
Department of Mathematics\\
Washington State University,  Pullman, WA 99164-3113 USA; \\
fax: (1 509) 335 1188; tel. (1 509) 335 3172; 
e-mail: khapala@math.wsu.edu

  \bigskip

\begin{abstract}
We study the global approximate controllability of the reaction-diffusion equation in a  parallelpiped  $ \Omega = (a_1,b_1 ) \times \ldots (a_n,b_n) \subset R^n $, governed by a multiplicative control in a reaction term.  It is assumed  that the initial state $ u_0 $ admits zeros only on the intersections of $ \Omega$ with finitely many hyperplanes, parallel to the sides of $ \Omega$, and that $ u_0$ changes   its sign after crossing such hyperplanes (we further refer to them as the ``hyperplanes of change of sign'' or ``zero hyperplanes''). 
This paper can be viewed as a continuation  of work  presented in \cite{CanKh, CanKh2} for the controllability of the one dimensional reaction-diffusion equation with  solutions admitting finitely many zeros. However, the methods of  \cite{CanKh, CanKh2} are intrinsically one dimensional, while in this paper we introduce  a  novel approach to deal with the case of multiple spatial variables.

\end{abstract}

\bigskip

{\bf 1. Introduction.} The traditional linear operator methods, based on the duality pairing, developed to study controllability of linear evolution systems with additive controls, do not apply to nonlinear control problems arising in the context of multiplicative controls. The latter also represent a principlly different class of applications (such as, e.g.,  chain-reactions,  see \cite{Kh3} and the references therein).

Among early works on  ``multiplicative'' controllability,  let us  mention a pioneering work  \cite{BMS} (1982) by Ball, Marsden and Slemrod, establishing the approximate controllability of the rod  and  wave equations, based on an  implicit nonharmonic Fourier series approach, adapted for the  time-dependent (only) multiplicative controls. 

In turn, the multiplicative  controllability of linear and semilinear parabolic equations in several spatial dimensions was originated in the series of papers by Khapalov \cite{Kh1}-\cite{Kh2},  further   summarized in monograph \cite{Kh3}. The approach of \cite{Kh3}   makes use of   explicit asymptotic qualitative  methods that employ  piecewise constant-in-time multiplicative controls of both spatial and time variables.

In more  recent papers \cite{CanKh}, \cite{CanKh2}  we discussed  the multiplicative  controllability properties of the one   dimensional reaction-diffusion equation with the initial and target states admitting finitely many ``matching''  changes of sign. The results were extended in \cite{CanKh2} to a two dimensional Heat equation {\em on a disk} in the case when its solutions are radially symmetrical.  Nonetheless, the methods of \cite{CanKh}, \cite{CanKh2}  deal with   simultaneous control of  the motion of isolated zero points and are {\em intrinsically one dimensional} (with respect to the spatial variable) as they focus on the controlling the sign of $ u_t$ at zero points.

In the case of several spatial dimensions, the sets of points of change of sign  of solutions to the pde at hand  are no  longer  finite sets of points. Such principal change of setup  demands  a new methodology, which is the subject of this paper. Namely, we intend to exploit the fact that the choice of multiplicative control affects the selection of eigenfunctions of the spectral problem at hand. Respectively, our novel strategy is based on dealing  with multiplicative controls for which the desirable target state becomes the first ``essentially'' non-zero eigenfunction in the Fourier  series expansion of the resulting solution.

\bigskip
 
{\bf Controlled evolution system.} In this paper we consider the following reaction-diffusion equation:
   $$
   u_t \; = \; \Delta u \; + \; v(x,t)  u  \;\;\;\;\;  {\rm in} \;\;\;
   Q_T = \Omega \times (0, T), \; \;\; T>0, \; x = (x_1, \ldots , x_n)\eqno(1.1)$$
   $$
   u \mid_{\partial \Omega}   = 0, \;\;\;\; t \in (0, T),  \;\;\;\;
   u\:\mid_{t = 0} \; = u_0 \in H^2(\Omega) \bigcap H^1_0 ( \Omega ),
    $$
where   $ \Omega = (a_1,b_1 ) \times \ldots (a_n,b_n) \subset R^n $. The symbol  $ v \in L^\infty (Q_T) $ stands for   the  {\em multiplicative control} function, which we further assume to be piecewise constant-in-time.

It is known that, for any $T>0$, system (1.1) admits a unique solution in \\
 $ H^{2,1} (Q_T) \bigcap C ([0, T]; H^2 (\Omega) \bigcap H_0^1 (\Omega))$. 
 Here and below, we use the standard notations for Sobolev spaces, in particular, $ \; H_0^1 (\Omega) = \{\phi \mid \; \phi, \phi_{x_i} \in L^2 (\Omega),  \; \phi \mid_{\partial \Omega}= 0, \; i = 1, \ldots, n\}$, and $ \; H^{2,1} (Q_T) = \{\phi \mid \; \phi, \phi_{x}, \phi_{x_i x_j }, \; \phi_t  \in L^2 (Q_T)\; i,j = 1, \ldots, n\},  \; H^2 (\Omega) = \{\phi \mid \; \phi, \phi_{x_i},\phi_{x_i x_j} \in L^2 (\Omega)\} \; i,j = 1, \ldots, n\}$.

Let us remind the reader that, in the classical sense, an evolution system is called {\em globally approximately
controllable in a given space $H$} at time
$ T > 0$, if   it can be steered in $ H$ from any initial state
into any neighborhood of any desirable target state  at time $ T$
by making use of  a suitable available  control. 
However,  this type of controllability is out of question  for  boundary problem (1.1) (e.g., 
the zero-function is the fixed point of the right-hand side in (1.1)).

The layout of the paper is as follows. In the next section we discuss some critical ``observations'' exposing the principal difficulties  which one can encounter in the framework of multiplicative controllability,  particularly,  due to the {\em  maximum principle} for solutions of parabolic equations.
In Section 3, we describe the proposed  methodology and state the main results of this paper. In Sections 4-7, we prove  the main results.

\bigskip

{\bf 2. Maximum principle and  multiplicative controllability of linear parabolic PDE's: some critical observations.}   
In the following observations and examples we assume that  $ \Omega = (0,1) \times (0,1)$
(though, of course, it is just for illustratuion purposes).

\bigskip

{\bf Observation 1:}   Assume that  in (1.1)  $ u_0 (x) \geq 0$. Then, the  classical maximum principle requires the respective solution to (1.1) to stay nonnegative  at any moment of time $ t> 0$, regardless of the choice of $ v$. This means that system (1.1) cannot be steered from any such $ u_0 $ to a target state which is negative on a nonzero measure set in  $\Omega$. However, one can hope to approximately steer it to any non-negative state   in $ L^2 (\Omega)$ and it was shown  in \cite{Kh1}-\cite{Kh3} for a rather  general class of {\em semilinear} parabolic equations in several spatial dimensions. 

\bigskip

{\bf Observation 2:} Let us assume now that $u_0$  admits   changes of sign across finitely many curves, say, $S_i(0), i = 1, \ldots, M$,  in  domain  $ \Omega$, splitting $ \Omega$ into  $N$ simply-connected open subdomains $ A_i(0), i = 1,\ldots, N$, in which  $u_0$ does not change sign. A straightforward adaptation of reasoning in Observation 1  implies that  for any control $ v$ the respective solution $ u$ to (1.1) at any moment of time $t>0$  cannot have {\em more than} as many curves of  sign change (while some curves of change of sign can ``merge''). Furthermore, the ``geometrical topology''  of these evolving curves  should be compatible with that of $ u_0$. Namely,
the respective $ A_i (t)$'s  are the result of ``continuous transformations'' of the original $ A_i (0)$'s, see the illustrating examples below on Figures 0-5 (the change of sign 
is indicated with ``$\pm$'' symbols).

\setlength{\unitlength}{1mm}
\begin{picture}(200,80)(0,0)
\linethickness{1pt}
  
\put (52,55){\vector(1,0){10}}

\put (92,55){\vector(1,0){10}}

\put(20,20.5){\framebox(30,30)}

\put(25,20.5){\line(0,1){30}}
\put(45,20.5){\line(0,1){30}}

\put(22,35){\makebox(0,0)[b]{$+$}}
\put(35,35){\makebox(0,0)[b]{$-$}}
\put(48,35){\makebox(0,0)[b]{$+$}}

\put(22,46){\makebox(0,0)[b]{$-$}}
\put(35,46){\makebox(0,0)[b]{$+$}}
\put(48,46){\makebox(0,0)[b]{$-$}}

\put(22,23){\makebox(0,0)[b]{$-$}}
\put(35,23){\makebox(0,0)[b]{$+$}}
\put(48,23){\makebox(0,0)[b]{$-$}}

\put(20,30){\line(1,0){30}}
\put(20,45){\line(1,0){30}}
\put(22,35){\makebox(0,0)[b]{$+$}}

\put(48,35){\makebox(0,0)[b]{$-$}}

\put(60,20.5){\framebox(30,30)}

\put(65,20.5){\line(0,1){30}}
\put(85,30){\line(0,1){21}}
\put(85,20.5){\line(0,1){6.5}}
\put(85,27){\line(1,0){4.8}}

\put(62,35){\makebox(0,0)[b]{$+$}}
\put(75,35){\makebox(0,0)[b]{$-$}}
\put(88,35){\makebox(0,0)[b]{$+$}}

\put(62,46){\makebox(0,0)[b]{$-$}}
\put(75,46){\makebox(0,0)[b]{$+$}}
\put(88,46){\makebox(0,0)[b]{$-$}}

\put(62,23){\makebox(0,0)[b]{$-$}}
\put(75,23){\makebox(0,0)[b]{$+$}}
\put(88,23){\makebox(0,0)[b]{$-$}}

\put(60,30){\line(1,0){25}}
\put(60,45){\line(1,0){30}}

\put(100,20.5){\framebox(30,30)}

\put(105,20.5){\line(0,1){30}}
\put(125,30){\line(0,1){20.5}}

\put(102,35){\makebox(0,0)[b]{$+$}}
\put(115,35){\makebox(0,0)[b]{$-$}}
\put(128,35){\makebox(0,0)[b]{$+$}}

\put(102,46){\makebox(0,0)[b]{$-$}}
\put(115,46){\makebox(0,0)[b]{$+$}}
\put(128,46){\makebox(0,0)[b]{$-$}}

\put(102,23){\makebox(0,0)[b]{$-$}}
\put(115,23){\makebox(0,0)[b]{$+$}}

\put(100,30){\line(1,0){25}}
\put(100,45){\line(1,0){30}}

\end{picture}

\setlength{\unitlength}{1mm}
\begin{picture}(200,80)(0,0)
\linethickness{1pt}

\put (12,55){\vector(1,0){10}}

\put(20,20.5){\framebox(30,30)}  

\put(25,30){\line(0,1){15}}
\put(45,30){\line(0,1){15}}

\put(22,35){\makebox(0,0)[b]{$+$}}
\put(35,35){\makebox(0,0)[b]{$-$}}
\put(48,35){\makebox(0,0)[b]{$+$}}

\put(35,46){\makebox(0,0)[b]{$+$}}

\put(35,23){\makebox(0,0)[b]{$+$}}

\put(22,35){\makebox(0,0)[b]{$+$}}

\put(48,35){\makebox(0,0)[b]{$-$}}

\put(25,30){\line(1,0){20}}
\put(25,45){\line(1,0){20}}

\put(75,8){\makebox(0,0)[b]{Figure 0.}}
\end{picture}

\bigskip

{\bf Observation 3:} Due to the maximum principle of solutions to parabolic equations, {\em no new zero-curves of change of sign  can emerge} inside any open  set where, prior to that, we had  $u  (x,t) \geq  0$ or $u  (x,t) \leq  0$ for all spatial points.

\bigskip

Let us illustrate the above ideas  with the following two `hypothetical''  examples.

\bigskip

{\sf Example 2.1.}  In (1.1) with $ \Omega = (0,1) \times (0,1)$ assume that $ u_0 (x)$ is a function, describing the initial temperature distribution $ u_0 (x)$,  is positive  on the left of a single vertical line and is negative  on the right of it as shown on the left square on Fig. 1.  The change of sign 
is indicated with ``$\pm$'' symbols. 

On the following two squares we can see two {\em possible}  positions of the zero line which  moves, as indicated with arrows,  to the right  towards the right vertical part of  the boundary $ \partial \Omega$ and finally it can ``merges'' with it (or  get arbitrarily  close to it).

\setlength{\unitlength}{1mm}
\begin{picture}(200,80)(0,0)
\linethickness{1pt}

\put (52,55){\vector(1,0){10}}

\put (92,55){\vector(1,0){10}}

\put(20,20.5){\framebox(30,30)}

\put(30,20.5){\line(0,1){30}}

\put(60,20.5){\framebox(30,30)}

\put(80,20.5){\line(0,1){30}}

\put(129,20.5){\line(0,1){30}}

\put(100,20.5){\framebox(30,30)}

\put(25,35){\makebox(0,0)[b]{$+$}}

\put(35,35){\makebox(0,0)[b]{$-$}}

\put(72,35){\makebox(0,0)[b]{$+$}}

\put(84,35){\makebox(0,0)[b]{$-$}}

\put(115,35){\makebox(0,0)[b]{$+$}}

\put(75,12){\makebox(0,0)[b]{Figure  1.}}

\end{picture}

\bigskip

{\sf Example 2.2.}  On Fig. 2 below   we again consider a squared  domain $ \Omega$ but now with two vertical and two horizontal zero-lines of change of sign for some initial temperature distribution  $ u_0 (x)$. The change of sign is again indicated with ``$\pm$''. 

On the three squares below we illustrate a possible evolution of positions of these zero lines  compatible with the above-mentioned maximum principle.

\setlength{\unitlength}{1mm}
\begin{picture}(200,80)(0,0)
\linethickness{1pt}
  
\put (52,55){\vector(1,0){10}}

\put (92,55){\vector(1,0){10}}

\put(20,20.5){\framebox(30,30)}

\put(25,20.5){\line(0,1){30}}
\put(45,20.5){\line(0,1){30}}

\put(22,35){\makebox(0,0)[b]{$+$}}
\put(35,35){\makebox(0,0)[b]{$-$}}
\put(48,35){\makebox(0,0)[b]{$+$}}

\put(22,46){\makebox(0,0)[b]{$-$}}
\put(35,46){\makebox(0,0)[b]{$+$}}
\put(48,46){\makebox(0,0)[b]{$-$}}

\put(22,23){\makebox(0,0)[b]{$-$}}
\put(35,23){\makebox(0,0)[b]{$+$}}
\put(48,23){\makebox(0,0)[b]{$-$}}

\put(20,30){\line(1,0){30}}
\put(20,45){\line(1,0){30}}
\put(22,35){\makebox(0,0)[b]{$+$}}

\put(48,35){\makebox(0,0)[b]{$-$}}

\put(60,20.5){\framebox(30,30)}

\put(65,20.5){\line(0,1){30}}
\put(85,20.5){\line(0,1){30}}

\put(62,38){\makebox(0,0)[b]{$+$}}
\put(75,38){\makebox(0,0)[b]{$-$}}
\put(88,38){\makebox(0,0)[b]{$+$}}

\put(62,46){\makebox(0,0)[b]{$-$}}
\put(75,46){\makebox(0,0)[b]{$+$}}
\put(88,46){\makebox(0,0)[b]{$-$}}

\put(62,27){\makebox(0,0)[b]{$-$}}
\put(75,27){\makebox(0,0)[b]{$+$}}
\put(88,27){\makebox(0,0)[b]{$-$}}

\put(60,35){\line(1,0){30}}
\put(60,45){\line(1,0){30}}

\put(100,20.5){\framebox(30,30)}

\put(115,20.5){\line(0,1){30}}
\put(125,20.5){\line(0,1){30}}

\put(102,38){\makebox(0,0)[b]{$+$}}
\put(120,38){\makebox(0,0)[b]{$-$}}
\put(128,38){\makebox(0,0)[b]{$+$}}

\put(102,46){\makebox(0,0)[b]{$-$}}
\put(120,46){\makebox(0,0)[b]{$+$}}
\put(128,46){\makebox(0,0)[b]{$-$}}

\put(102,26){\makebox(0,0)[b]{$-$}}
\put(120,26){\makebox(0,0)[b]{$+$}}

\put(100,35){\line(1,0){30}}
\put(100,45){\line(1,0){30}}

\put(75,12){\makebox(0,0)[b]{Figure  2.}}

\end{picture}

\bigskip
 {\bf 3. Methodology  and  main results.} 
In this paper we will introduce a novel approach to multiplicative controllability  exploiting the  idea  to employ  multiplicative controls which would make the desirable target states to be co-linear to  the first ``essentially non-zero'' eigenfunction in the Fourier  series expansion of the resulting solutions.

This new strategy will become possible due to the ``supporting'' qualitative techniques, introduced in  \cite{Kh1}, \cite{Kh2} (see also \cite{Kh3}) in the context of multiplicative controllability,  exploiting     the idea  of extracting various  qualitative dynamics from system (1.1),  namely, by making use of: 
\begin{itemize}
\item
either ``long-term'' static controls (i.e., not depending on the time-variable), which will drive system (1.1) to the first non-zero term in the Fourier series expansion of the respective solution (while the  subsequent terms are dissipating in time), 
\item
or  static controls acting on ``vanishingly small''  time-intervals to ensure that  {\em the reaction dynamics} in (1.1) will dominate over the diffusion process, when such controls are applied. 
\end{itemize}

In our main results below we will use the following  definition detailing the concept of change of sign of a function in a domain $\Omega \subset R^n, n = 2, \ldots$  that of our interest in this paper.

\bigskip
{\bf Definition 3.1.}  {\em Everywhere in this paper, when we talk  about ``finitely many $(n-1)$-dimensional surfaces of change of sign'' or ``zero-surfaces/hyperplanes'' for a given  function,  we mean that the total measure of  these  surfaces is zero and that they can split $ \Omega$ into  finitely many subdomains (simply-connected open sets) $ S_i^+, i = 1, \ldots, l  $ and $ S_j^-, j = 1, \ldots, m $, bounded by the respective aforementioned surfaces and, possibly,  by parts of boundary $ \partial \Omega$. In addition, they are  such that the function at hand  is positive inside of any of  $ S_i^+, i = 1, \ldots, l  $  and is negative  inside  of any of $ S_i^-, i = 1, \ldots, m $  (except, possibly, on sets of zero measure). We also assume that  this function must change its sign (in the above sense, that is, almost everywhere) across these surfaces.}

\bigskip

{\bf Main methdological question.} In the previous section we pointed out at some principal  restrictions  on the target states which can be reached from a  \underline{given initial state} $u_0$ along the dynamics of (1.1).  Respectively, in this paper, we intend to attack this problem from the ``opposite direction'', namely, we  ask:  

\bigskip
\par\noindent {\em \underline{From what initial states $u_0 $ system (1.1) can  be steered to a given desirable target  state} $u_1$?}
 
  \bigskip
Our first result deals with a special  case when curves  of change of sign do not need to be moved.

 {\bf Theorem 3.1.}  {\em  Assume that the initial state $u_0 \in H^2 (\Omega) \bigcap H^1_0 (\Omega$ in (1.1) and the desirable target states $ u_1 \in L^2 (\Omega)$  have the same surfaces   of change of sign (see Definition 3.1). Then system (1.1) can be steered from  the former   to  the latter  as close as we wish in  $ L^2 (\Omega)$ at some time $ T>0$.}

\bigskip
Theorem 3.1  is proven in Appendix. It is  instrumental for the proofs of Theorems 3.3-3.6, particularly, due to the following remark.

\bigskip
{\bf Remark 3.1.} {\em In Theorem 3.1 the steering is achieved as $ T \rightarrow 0+$}.
 
\bigskip

{\bf Remark 3.2.}  {\em We assume that solutions of all introduced below spectral problems are orthonormalized in the respective $L^2$-spaces.}

\bigskip

Let  $ \{\lambda_k, \omega_k\}_{k = 1}^\infty $ denote the eigenelements of the spectral problem associated with (1.1):
$$
\Delta \omega_k  + v_0(x) \omega_k = \lambda_k \omega_k, \;\; \lambda_1 \geq  \lambda_2 \geq \ldots, \;\;\;\; \lambda_k \rightarrow -\infty \;\; {\rm as} \;  k \rightarrow -\infty.
\eqno(3.1)$$
Then, the solution to (1.1) with $ v = v_0  $ admits the following Fourier series representation:
$$
u (x, t) \; = \;\; \sum_{k = 1}^{\infty}  \left( \int_\Omega u_0 \omega_k dx \right) e^{\lambda_k t} \omega_k(x).
\eqno(3.2)$$

For each  positive integer $ k_*$  introduce the following ``$k_*$-momentum problem''.

\bigskip
{\bf $k_*$- Momentum problem.} {\em  Let  in (3.1) $ \lambda_{k_*} \; >  \; \lambda_{k_*+1}$. Find  $ u_0^*$   such that:
$$
 \int_\Omega u_0^*  \omega_k dx = 0, \;\; k = 1, \ldots, k_*-1,
 \eqno(3.3{\rm a})$$
 $$
 \int_\Omega u_0^* \omega_{k_*}  dx \; =  \; c_0 \neq 0  
 \eqno(3.3{\rm b})$$
for some $ c_0$. 
}

\bigskip

{\bf Remark 3.3.} {\em We would like to emphasize here the strict nature of the inequality $ \lambda_{k_*} \; >  \; \lambda_{k_*+1}$.
}

\bigskip

The next result, Theorem 3.2, 
 describes one of the central ideas of  our  controllability strategy in this paper.

\bigskip

{\bf Theorem 3.2.}  {\em Consider any $ v_0 \in L^\infty (\Omega)$ and a positive integer $ k_* > 1$. Let, for this $ v_0$  $ \lambda_{k_*} > \lambda_{k_*+1} $ in (3.1). Assume that  the $k_*$- momentum problem (3.3a-b) admits a solution  for some $ u_0 \in H^2 (\Omega) \bigcap H^1_0 (\Omega)$, that is,
$$
u_0 = \sum_{k= k_*}^\infty  a_k \omega_k, \; a_{k*} \neq 0.
$$ 
Then any of the target states $ \alpha  \omega_{k_*} , \alpha >0 $ can be approximately reached in $ L^2 (\Omega)$ from this $ u_0$  at some time $ T \geq T_1$ along the dynamics of  (1.1) with $ v = v_0 (x)$. }

\bigskip

{\bf Remark 3.4} {\em If $ k_* = 1$, then the first eigenfunction is $ \omega_{k_*}$ is not multiple (due to the aforementioned classical maximum principle for solutions of  parabolic equations) and it does not change its sign in $ \Omega$, say, is non-negative). Hence, the respective $k_*=1$- momentum problem (3.4a-b) becomes  trivial.  In this case, with the help of Theorem 3.1, the result of Theorem 3.2 holds for any non-negative target state and, thus,  we have  the ``non-negative'' approximate controllability (to all non-negative target states) in $ L^2 (\Omega)$ as it was shown earlier in  \cite{Kh1}-\cite{Kh2}}.

\bigskip
{\bf Remark 3.5.} {\em Theorems 3.1 and 3.2 do not use the specific geometric structure of $\Omega$ and, thus, apply to general domains.}


\bigskip
Let us assume that our initial state $ u_0 (x), x =(x_1, \ldots, x_n) $  changes  its sign in $ \Omega$  on $k_i-1 \geq 1 $ hyperplanes $ Q_{ij}$ perpendicular to the $x_i-$axis
 $$
  Q_{ij} = \{ x =(x_1,  \ldots , x_n) | \; x_i =  x^0_{ij} \},  \; a_i <  x^0_{i1} < \ldots < x^0_{i k_i-1} < b_i, \;    i = 1, \dots, n, j =  1, \ldots, k_i -1.
\eqno(3.4{\rm a}) $$

Select any set of   functions  $ u_{i} (x_i) \in L^2 (a_i, b_i), i = 1, \ldots, n$ such that:
\begin{itemize}
\item
each $ u_{i} (x_i) $ has exactly $ k_i-1$ points of change of sign at $ a_i <  x_{i1} < \ldots < x^0_{i k_i-1} < b_i$ 
\item
and the function
$$
u^* (x) = \Pi_{i=1}^n   u_{i} (x_i)
\eqno(3.5{\rm a})$$
has the same sequence of change sign as $ u_0$.
\end{itemize}
 In the above and below, if $ k_i = 1$, then there is no zero-hyperplane perpendicular the the $i$-dimension.
 
In turn, assume that our target state $ u_1 (x), x =(x_1, \ldots, x_n) $ changes  its sign in $ \Omega$  on $k_i-1  \geq   0 $ hyperplanes $ P_{ij}$ (or ``zero-hyperplanes") perpendicular to the $x_i-$axis
 $$
  P_{ij} = \{ x =(x_1,  \ldots , x_n) | \; x_i =  x_{ij} \},  \; a_i <  x_{i1} < \ldots < x_{i k_i-1} < b_i, \;    i = 1, \dots, n, j = 1, \ldots, k_i -1.
\eqno(3.4{\rm b}) $$

Select a set of   functions  $ w_{i} (x_i), i = 1, \ldots, n $ such that they are either constant (say, identically equal to 1) or:
\begin{itemize}
\item
$$
w_{i} (x_i)  \in H^2 (a_i,b_i) \bigcap H_0^1 (a_i,b_i) \bigcap C^2[a_i, b_i], \;    i = 1, \dots, n;
$$ 
\item
each $ w_{i} (x_i) $ has exactly $ k_i-1$ points of change of sign at $ a_i <  x_{i1} < \ldots < x_{i k_i-1} < b_i$;
\item
selection of $w_i(x_i) \equiv 1$ takes place if there are no zero-hyperplanes perpendicular to the $i$-th axis,
\item
and the function
$$
w(x) = \Pi_{i=1}^n   w_{i} (x_i)
\eqno(3.5{\rm b})$$
has the same sequence of change sign as $ u_1$ in (3.4b);
\item
each $w_i$ is linear near $x_{i j}$'s (whence, $w_{ix_ix_i} =0$ near these points). 
\end{itemize}

Select control $ v=v_0$ for (1.1) as follows:
$$
v_0 (x)  \; = \; \sum_{i=1}^n v_i (x_i) , \;\;
v_{i} (x_i) = -\frac{w_{ix_i x_i} (x_i)}{w_{ix_i}} \in C [0, 1], 
\eqno(3.6)$$
at finitely many points (see Remark 5.2) where the denominator is not zero, and set $ v_i = 0$ otherwise.

Denote by  $\omega_{ij}$'s  and $\lambda_{ij}$'s the (orthonormalized)  solutions of  the following  spectral problem in $ H^2 (a_i, b_i) \bigcap H^1_0 (a_i, b_i) \bigcap C^2[a_i, b_i]$:
$$
 \omega_{ij x_i x_i } (x_i) + v_i  (x_i)  \omega_{ij}(x) = \lambda_{ij} \omega_{ij} (x_i), \; i = 1, \ldots, n, \; j = 1, \ldots
\eqno(3.7{\rm a})$$

{\bf Remark 3.6: Instrumental  observations.} 
\begin{itemize}
\item
{\em All ``one-dimensional'' eigenvalues $\lambda_{ij}$'s are \underline{simple}: }
$$
\lambda_{i1} <  \lambda_{i2} < \ldots,  \lambda_{ij}  \rightarrow \infty \;\; {\rm as} \;  j \rightarrow -\infty, \; i = 1, \ldots, n.
$$
\item
{\em  \underline{ Note that $w_i (x_i)$ is co-linear to $ \omega_{i k_i} (x_i)$ and $\lambda_{ik_i} = 0$} (see (3.6)). Indeed, $ \omega_{i k_i} (x_i)$ has   $ k_i-1$ zeros in $ (a_i, b_i)$, $ \omega_{i k_i-1 } (x_i)$ has $ k_i-2$ zeros and so forth ... (see Remark 5.2 below).}
\end{itemize}

\bigskip
Theorems 3.1 and 3.2 do not use the specific geometric structure of $\Omega$ and, thus, apply to general domains. However, our next result fully exploits the fact that $\Omega$ is a parallelepiped in order to apply the method of separation of variables to reduce the multidimensional moment problem in (3.3a-b) to a set of one dimensional moment problems. It also combines the results of Theorems 3.1 and 3.2.

\bigskip
{\bf Theorem 3.3.} {\em Let $u_0 \in H^2(\Omega) \bigcap  H^1_0 (\Omega)$ have  zero-hyperplanes (of change of sign) as in (3.4a) and $u_1 \in L^2 (\Omega)$ have   zero-hyperplanes as in (3.4b). Let    $w_i(x_i), i = 1, \ldots, n$ be any set of functions as in (3.5b) (there are infinitely many sets like that). Assume that there exist functions  $u_i, i = 1, \ldots, n$ from the closure in $ L^2 (\Omega)$ of the set of functions as in (3.5a) such that:
\begin{itemize}
\item
$$
 \int_{a_i}^{b_i} u_{i} (x_i)  \omega_{ij} (x_i) dx_i  \; = \; 0, 
 \; j= 1, \ldots, k_{i} -1, \; i = 1, \ldots, n,
 \eqno(3.7{\rm b}) $$
 or, alternatively, we can fund   sequences $\{u^s_i \in L^2 (a_i, b_i)\}_{s=1}^\infty, i = 1, \ldots, n$  as in (3.5a) such that 
 $$
 \lim_{s \rightarrow \infty} \int_{a_i}^{b_i} u^s_{i} (x_i)  \omega_{ij} (x_i) dx_i  \; = \; 0, 
 \; j= 1, \ldots, k_{i} -1, \; i = 1, \ldots, n,
 \eqno(3.7{\rm c})$$ 
 \item
 while
 $$
 \Pi_{i=1, \ldots, n}    |  \int_{a_i}^{b_i} u_{i} (x_i)  \omega_{i k_i} (x_i) dx_i  | \; = \; 1,
\eqno(3.8{\rm a})$$
or, respectively, 
 $$
 \Pi_{i=1, \ldots, n}    |  \int_{a-i}^{b_i} u^s_{i} (x_i)  \omega_{i k_i} (x_i) dx_i  | \; = 1,  \;\;{\rm as} \; s \rightarrow \infty,
\eqno(3.8{\rm b})$$
independently of the choice of the ``degree of closedness to zero'' of terms in (3.7c). 
\end{itemize}
Then, system (1.1)  can  be steered from $ u_0$ to $ u_1$ as close as we wish in $ L^2 (\Omega)$ at some time $ T>0$.}

\bigskip
{\bf Discussion of conditions (3.7b)-(3.8b).}  These conditions form a set of moment problems in one spatial dimension for each $i = 1, \ldots $ on a cone (that is, these are  {\bf not}  standard linear moment problems  in a Hilbert space), derived from the moment problem  (3.3a-b). Namely, for each dimension $i$ we want to find a sequence of functions $\{u_i^s\}_{s=1}^\infty$ which:
\begin{enumerate}
\item
in general, ``approaches''  a subspace of $L^2 (a_i, b_i)$ perpendicular to the span of  $\{\omega_{ij}\}_{j = 1}^{k_{i*}-1}$ (it may {\bf not} be a converging sequence), namely, containing functions  of the form 
$$
\sum_{k= k_{i*}}^\infty  a_k \omega_{ij} (x_i).
$$ 
In other words, $u_i^s$'s tend to become ``perpendicular to the aforementioned  span (see (3.7b-c));
\item
but {\bf not} to become perpendicular to $\omega_{ik*}$ (see  (3.8a-b));
\item
functions $u_i^s$'s  should {\em all have the same zeros in the dimension $i$ as the given initial condition $ u_0$ in (1.1) with the same sequence of change of sign}.
\item
If (3.7b-c)-(3.8a-b) hold, then the first ``essentionaly non-zero'' terms in the Fourier  series representations of respective solutions to (1.1), (3.6), with any of the initial conditions described in Therem 3.3, have the same sequence of zero-hyperplanes as any of the targets in this theorem, that is, as $ u_1$.
\end{enumerate}

\bigskip

We will show below that Theorem 3.3, in particular,  implies  following straightforward implication in the case when (3.7b-c)-(3.8a-b) become trivial.

\bigskip

\bigskip
{\bf Theorem 3.4.} {\em System (1.1) can  be steered from  any  initial state $ u_0 \in H^2(\Omega) \bigcap H^1_0 (\Omega)$, with {\bf at most one, per each dimension $i$,}  hyperplane of change of sign as in (3.4a), to any target state $u_1 \in H^2(\Omega) \bigcap  H^1_0 (\Omega)$, which has a respectively equal amount of  hyperplanes of change of sign, that is,  as in (3.4b), as close as we wish in $ L^2 (\Omega)$ at some time $ T>0$. }

\bigskip

In the general case of finitely many hyperplanes of change of sign we have the following two results.

\bigskip
{\bf Theorem 3.5.} {\em Consider any target state $u_1 \in H^2(\Omega) \bigcap  H^1_0 (\Omega)$ with given finitely many hyperplanes of change of sign (as in (3.4b)). Then system (1.1)  can  be steered to $ u_1$, as close as we wish in $ L^2 (\Omega)$,  from  any  initial state $ u_0 \in H^2(\Omega) \bigcap H^1_0 (\Omega) $, which has respectively  the same amount of similarly oriented  hyperplanes of change of sign (as in (3.4a)), provided that the points $ x^0_{ij}$'s in (3.4a) are selected arbitrarily in their respective segments $(a_i, b_i)$, except for, possibly, a set of at most countably many points.}

\bigskip

The proof of Theorem 3.5 deals with the following  two assumptions.

\bigskip

{\bf Assumption 3.1.}  {\em The following sets of vectors are linear independent in the respective space $R^{k_i-1}$:
$$
\{(\omega_{ij} (x^0_{i1}), \ldots, \omega_{ij} (x^0_{ik_i-1}))\}, \;  i = 1, \dots, n, \; j = 1, \ldots, k_i -1,
\eqno(3.9)$$
where $x^0_{ij}$ are from (3.4a).}

\bigskip

\bigskip

{\bf Assumption 3.2.}  {\em If Assumption 3.1 does not hold for some (dimension) $i$, that is, the respective set in (3.9) is linear dependent, assume that
$$
\{(\omega_{ik_i} (x^0_{i1}), \ldots, \omega_{ik_i -1} (x^0_{ik_i}))\} \not\in {\rm span} \\\, \{\{(\omega_{ij} (x^0_{i1}), \ldots, \omega_{ij} (x^0_{ik_i -1}))\} |   \;  j = 1, \ldots, k_i -1\}.
\eqno(3.10)$$
}

\bigskip
{\bf Theorem 3.6.} {\em The statements of  Theorem 3.5 hold under Assumptions 3.1 and 3.2.}

\bigskip
{\bf Discussion of Assumptions 3.1.} Note  that in case of Theorem 3.4, Assumption 3.1 holds in a trivial way for $ k_i -1  = 1$  - we just have one non-zero vector  in (3.9).  

 Let us assume that Assumption 3.1 does not hold for some dimension $i$. We claim that the set of $x_i$'s for which it fails cannot have interior points in $ (a_i, b_i)$. 
To this end, we will investigate the linear dependence of vectors forming the columns of matrix generated my the vectors in (3.9).

Let, for some $ \alpha_j, j = 1, \ldots, k_i -1$, in some open interval $(c,d) \subset (a_i, b_i)$ he have: 
$$
g(x_i) = \sum_{ j = 1}^{k_i -1} \alpha_j \omega_{ij} (x_i ) \equiv 0.
\eqno(3.11)$$ 
Apply the elliptic operator in (3.7a) to (3.11)   $ \; k_i -1$ times to obtain the following linear algebraic system in $ \alpha_j \omega_{ij} (x_i ), j = 1, \ldots, k_i -1$ with Vandermonde matrix:
$$
\sum_{ j = 0}^{k_i -1} \alpha_j \omega_{ij} (x_i ) \equiv 0, \ldots, \sum_{ j = 0}^{k_i -1} \lambda^{k_i -1}_{ij} \alpha_j \omega_{ij} (x_i ) \equiv 0
$$
Since all $\lambda_{ij}$'s are single, it has only the trivial solution in $(c,d)$, that is, 
$$ 
\alpha_j \omega_{ij} (x_i ), j \equiv 0, \ldots, k_i -1, \; x_i \in (c,d).
$$
However,  $\omega_{ij} (x_i )$'s can have only finitely many zero's in $ (a_i, b_i)$  (see  Remark 5.1). Hence, all
$\alpha_j = 0$. Contradiction.

Thus, the set of $x_i$'s for which Assumption 3.1 fails, at most,  consists of points that are separated by open segments of points for which Assumption 3.1 holds.  In other words, the set of $x_i$'s, for which Assumption 3.1 fails, is {\em at most countable}.

\bigskip

Furthermore, if the functions $\omega_{ij} (z )$'s are analytic in $[a_i, b_i]$, then the set of $x_i$'s, for which Assumption 3.1 fails, is {\em finite}. Indeed, otherwise, the analytic function $ g(x_i)$ in (3.11) would be vanishing on a set of real points that have a limit point in $[a_i, b_i]$. Hence, $g(x) \equiv 0$ in $[a_i, b_i]$, which contradicts to Remark 5.1.

\bigskip
{\bf Example to Theorems 3.5-3.6: The case of two zero-hyperplanes.}  Let the initial state $u_0$ and target state $ u_1$ have two hyperplanes of change of sign perpendicular to  the $x_1$-axis, as described, respectively, in (3.4a) and (3.4b). Then, according to Remark 5.1, we have the following layout of the zero-points po change of sign for $ \omega_{11}, \omega_{12}$ and $ \omega_{13}$ (in the notattions of (3.4b)):
$$
a_1 < x_{11} < x_1^{12} < x_{12} < b_1,
$$
where $x_{11}$ and $x_{12}$ are zeros of $ u_1(x_1)$ and $\omega_{13} (x_1)$, $ x_1^{12}$ is the only zero-point of change of sign of $\omega_{12}(x_1)$, while $ \omega_{11}(x_1)$ does not have such zero-points. Without loss of generality, we can assume that immediately on the right of $ a_1$, the above eigenfunctions are positive.

Then, Assumption 3.1 holds if, e.g., the zeros of $u_0$ in (3.4a) are located as follows: 
$$ 
x^0_{11} \in (a_i, x_1^{12}), \;\; x^0_{12} \in (x_1^{12}, b_1)).
$$ 
Indeed, in this case components of vector $ (\omega_{11} (x_{11}^0), \omega_{11} (x_{12}^0))$ are of the same sign, while components of vector $ (\omega_{12} (x_{11}^0), \omega_{12} (x_{12}^0))$ have the opposite signs.

Alternatively, Assumption 3.1 can (theoretically)  fail  in the areas where  two aforementioned vectors have both coordinates of the same sign, {\em while  these vectors  have to be  co-linear}, say, 
$$ 
x^0_{11}, x^0_{12} \in (a_1, x_1^{12}).
$$ 
In this case, if 
$$ 
x^0_{11} \in (a_i, x_{11}), \;\; x^0_{12} \in (x_{11}, x_1^{12}),
$$
then vector $ (\omega_{13} (x_{11}^0), \omega_{13} (x_{12}^0))$  has  coordinates of the opposite sign and condition (3.10) hods, i.e., Assumption 3.2 holds.

\bigskip

{\bf Remark 3.7: Selection of zero-hyperplanes for $ u_0$ in Theorem 3.5.} {\em From the above discussion it follows that, for each dimension $i$, we can select the zero points  $ x^0_{ij}, j =1, \ldots, k_i-1$ for $u_o$ in Theorem 3.5 as follows;
\begin{itemize}
\item
Select $ x^0_{i1}$ on (3.4a) arbitrarily in $(a_i, b_i)$.;
\item
in order to make vectors
$$
\{(\omega_{ij} (x^0_{ij}), \ldots, \omega_{ij} (x^0_{ij}))\}, \;  j = 1, 2, 
$$
linear independent we can select the point  $ x^0_{i2}$ everywhere in $(a_i, b_i)$ except of a set of at most countably many points;
\item
and so on ...
\end{itemize}
}

\bigskip

{\bf Remark 3.8: About extension of Theorem 3.5.} {\em At a first glance, it may seem that one can apply some ``density argument'' to get rid of Assumption 3.1 in Theorems 3.5 and 3.6. However, if we select a sequence of auxiliary $u_{0s}$, satisfying Assumption 3.1, to approximate the actual $ u_0$, not satisfying this assumption, we may have a divergent sequence of controls associated with each $u_{0s}$, that is, in the framework of our methods used to prove Theorem 3.6.  These methods determine a suitable sequence of  controls  as solutions to a suitable linear algebraic system, with matrices constructed out of vectors in (3.9), i.e.,  with  non-zero determinants, guaranteed by Assumption 3.1.  In a ``conventional density argument'', if zeros of $u_{0s}$ approximate those of $u_0$, these determinants may converge to a degenerate one, associated with zeros of $u_0$ in (3.9).}

\bigskip

{\bf Remark 3.9: Geomentry of sero surfaces during the steering.} {\em  For the one dimensional case, the method  of his paper can be viewed as an alternative method to proof the main results in \cite{CanKh}-\cite{CanKh2} under the Assumptions 3.1 and 3.2. These assumptions  are not required in the aforementioned works. The main results in \cite{CanKh}-\cite{CanKh2} were achieved by finitely many (continuous) incremental moves of zero points.  In particular, this fact  be used to move zero points into open sets for which Assumption 3.1 holds, which would allow one to assume Assumption 3.1 without loss of generality. In the multidimensional case, considered in this paper, an analogous approach would be to move the zero-hyperplanes, which would also need to somehow maintain the  strict geometry of these hyperplanes at every moment of time. In the arguments of this paper we only require the zero surfaces to be hyperplanes perpendicular to the respective axes at the initial and final moments of steering.}


\bigskip
 
The following two figures illustrate Theorems 3.5 and 3.6  for the case of $\Omega = (0, 1) \times (0,1)$.

\setlength{\unitlength}{1mm}
\begin{picture}(200,80)(0,0)
\linethickness{1pt}

\put (52,55){\vector(1,0){10}}


\put(20,20.5){\framebox(30,30)}

\put(30,20.5){\line(0,1){30}}

\put(35,20.5){\line(0,1){30}}

\put(60,20.5){\framebox(30,30)}

\put(80,20.5){\line(0,1){30}}

\put(85,20.5){\line(0,1){30}}

\put(60,20.5){\framebox(30,30)}


\put(25,35){\makebox(0,0)[b]{$+$}}

\put(32,35){\makebox(0,0)[b]{$-$}}

\put(40,35){\makebox(0,0)[b]{$+$}}

\put(72,35){\makebox(0,0)[b]{$+$}}

\put(82,35){\makebox(0,0)[b]{$-$}}

\put(88,35){\makebox(0,0)[b]{$+$}}

\put(115,35){\makebox(0,0)[b]{$+$}}

\put(35,15){\makebox(0,0)[b]{$u_0$}}
\put(75,15){\makebox(0,0)[b]{$u_1$}}

\put(75,8){\makebox(0,0)[b]{Figure  3.}}

\end{picture}

\setlength{\unitlength}{1mm}
\begin{picture}(200,80)(0,0)
\linethickness{1pt}

\put (52,55){\vector(1,0){10}}


\put(20,20.5){\framebox(30,30)}

\put(30,20.5){\line(0,1){30}}

\put(20,30){\line(1,0){30}}

\put(60,20.5){\framebox(30,30)}

\put(80,20.5){\line(0,1){30}}

\put(60,40){\line(1,0){30}}

\put(60,20.5){\framebox(30,30)}


\put(25,35){\makebox(0,0)[b]{$+$}}

\put(40,35){\makebox(0,0)[b]{$-$}}

\put(25,25){\makebox(0,0)[b]{$-$}}

\put(40,25){\makebox(0,0)[b]{$+$}}


\put(72,45){\makebox(0,0)[b]{$+$}}
\put(72,30){\makebox(0,0)[b]{$-$}}

\put(84,45){\makebox(0,0)[b]{$-$}}

\put(84,30){\makebox(0,0)[b]{$+$}}


\put(35,15){\makebox(0,0)[b]{$u_0$}}
\put(75,15){\makebox(0,0)[b]{$u_1$}}

\put(75,8){\makebox(0,0)[b]{Figure  4.}}

\end{picture}

\bigskip

{\bf 4. Proof of Theorem 3.1.}  Consider any target state $u_1 = \alpha \omega_{k_*}$, where {\em $\alpha > 0$ is fixed}.

Solution to (1.1) with $ v = v_0 (x) $ and  $ u_0$ solving the $ k_*$-momentum problem admits the following Fourier series representation as in (3.2):
$$
u (x, t) \; = \; \sum_{k = 1}^{k* -1}  \left( \int_\Omega u_0 \omega_k dx \right) e^{\lambda_k t} \omega_k(x) 
$$
$$
+ \; \left( \int_\Omega u_0 \omega_{k*} dx \right) e^{\lambda_{k*} t} \omega_{k*} (x) \; 
+ \; \sum_{k = k* +1}^\infty  \left( \int_\Omega u_0 \omega_k dx \right) e^{\lambda_k t} \omega_k(x) 
$$
$$
=  \; \left( \int_\Omega u_0 \omega_{k*} dx \right) e^{\lambda_{k*} t} \omega_{k*} (x)  \; 
+ \; \sum_{k = k* +1}^\infty  \left( \int_\Omega u_0 \omega_k dx \right) e^{\lambda_k t} \omega_k(x)..
\eqno(4.1)$$

In the definition of $ k_*$-momentum problem we assumed that $\lambda_{k*} > \lambda_{k* + 1}$. 
Therefore, we can select a  control of the following form:
$$ 
v (x) = v_0 (x)  - \lambda_{k*} + a,
$$
where $a >0$ will be selected below in (4.4).

With this new control, the solution to (1.1) will take the form similar to (4.1), with the same eigenfunctions and the eigenvalues shifted by $ \lambda_{k*} -a$ units to the right:
$$
u (x, t) \; = \;  \; e^{at} \left( \int_\Omega u_0^* \omega_{k*} dx \right)  \omega_{k*} (x) \; 
+ \; e^{at} \sum_{k = k* +1}^\infty  \left( \int_\Omega u_0^* \omega_k dx \right) e^{(\lambda_k -\lambda_{k*})t} \omega_k(x).
\eqno(4.2)$$
Without loss of generality we can ssume that 
 $$
0 <  \alpha \; \neq \;  \int_\Omega u_0 \omega_{k*}  dx  \;  \; =  \; c_0  > 0,  
 \eqno(4.3)$$
see (3.3b) for $ c_0$. (Alternatively, we would select $ a = 0$ in (4.2)).

Select an arbitrary  sequence of  $ 0< T_k \rightarrow \infty$ as $ k \rightarrow  \infty$ and set 
$$
 a = a_i = \frac{1}{T_i} \ln (\frac{\alpha}{\int_\Omega u_0 \omega_{k*}  dx}) \; \rightarrow  0+ \;\; {\rm as} \; i \rightarrow \infty 
\eqno(4.4) $$
in which case:
$$
e^{aT_i} = \frac{\alpha}{\int_\Omega u_0 \omega_{k*}  dx}.
$$
In view of (4.3), the respective sequence of solutions to (1.1) at times $T_i$'s  will converge to $ u_1 = \alpha \omega_{k*}$ in $ L^2 (\Omega)$ as $ i \rightarrow \infty$. 

This ends  the proof of Theorem 3.2.

\bigskip
{\bf Remark 4.1.} {\em Note that, in the above proof,  the result of Theorem 3.2 is achieved as the time of steering tends to zero.}

\bigskip

{\bf 5. Proof of Theorems 3.3 and  3.4: The case of $ \Omega \subset R^2 $  and a single line of change of sign.} Without loss of generality (and for the sake of simplicity of notations), we can assume that
$$
\Omega = (0, 1) \times (0, 1).
$$
Let $u_0$ have just  one vertical line of change of sign. 
We intend to  show how one  can steer  system (1.1) in $L^2 (\Omega) $  from any such $u_0$ to any  state $ u_1$, which also has a single vertical line of change of sign at any desirable position within $ \Omega$, e.g., as shown, e.g.,   on Fig. 5:

\setlength{\unitlength}{1mm}
\begin{picture}(200,80)(0,0)
\linethickness{1pt}

\put (52,55){\vector(1,0){10}}

\put(20,20.5){\framebox(30,30)}

\put(30,20.5){\line(0,1){30}}

\put(60,20.5){\framebox(30,30)}

\put(85,20.5){\line(0,1){30}}

\put(25,35){\makebox(0,0)[b]{$+$}}

\put(40,35){\makebox(0,0)[b]{$-$}}

\put(73,35){\makebox(0,0)[b]{$+$}}
\put(35,15){\makebox(0,0)[b]{$u_0$}}
\put(75,15){\makebox(0,0)[b]{$u_1$}}
\put(88,35){\makebox(0,0)[b]{$-$}}

\put(75,8){\makebox(0,0)[b]{Figure  5.}}

\end{picture}

Let $ u_1$ have a sign change on a vertical line positioned at  $ x_1 = x_{11}$ ($ x = (x_1, x_2)$), while $ u_0$ changes its sign on a vertical line positioned at $ x_1 = x^0_{11} \neq x_{11}$ (otherwise, we can apply Theorem 3.1).
The plan of proof is as follows:

\begin{itemize}
\item
In {\bf Step 1} we select a stationary control $ v(x) = v_0(x_1) $ such that  one of the eigenfunctions $ \omega_{k_*} $ of  the respective spectral problem as in  (3.1)    has a single vertical line of sign  change at $x_1 = x_{11}$, the same as  the desirable target state $ u_1 (x)$.
\item
In {\bf Step 2}, making use of Theorem 3.1,  we will steer (1.1) on some $ (0, T_*)$, where $ T_*$ can be selected as small as we wish, to some intermediate auxiliary target state $u_*(x)$, which can ``approximately'' solve the $k_*$-momentum problem (3.2)-(3.3a-b)  for  the eigenfunctions associated with $ v(x) =  v_0(x_1) $. Thus, $u_*$ will have the same zero-line as $u_0$.
\end{itemize}

The aforementioned steering in Step 2 will ensure that, for $ t >  T_*$, the term containing $ \omega_{k_*} $  will be the first ``substantially nonzero, dominating'' term in the respective expansion of solution to system (1.1),  {\em if control $v_0 (x_1)$  from Step 1 is engaged}.

\begin{itemize}\item
In {\bf Step 3}, making use of the argument of Theorem 3.2, we will show that we can  steer system (1.1)  from $ u_*$  to a state co-linear with $ \omega_{k_*} $ as close as we wish in $ L^2 (\Omega)$. 
\item
In {\bf  Step 4}  we will apply Theorem 3.1 again to further (approximately in $ L^2 (\Omega)$) steer system (1.1)   to the target state $u_1(x)$ without change of the position of the zero-line at $ x_{11}$.
\end{itemize}

\bigskip

{\bf Step 1: Selection of an auxiliary target state $ w (x_1) $ and associated control $ v = v_0(x)$.} 

Select, in notations of (3.5b),  any  function $ w_1  \in H^2 (0,1) \bigcap H_0^1 (0,1) \bigcap C^2[0,1]$,
$w_2 (x_2) \equiv 1$,  (see (3.5b)) such that 
$$
w_1  (x_1) = 
\left\{ \begin{array}{ll}
0, \;\;\;\; & \; {\rm for} \;  x = 0, x_{11}, 1,   \\ 
\neq 0, \;\;\;\; & \; {\rm for} \;  x \neq 0, x_{11}, 1,
\end{array}
\right.
|| w||_{L^2 (0,1)} = 1, 
$$
whence:
$$
w(x) = w_1 (x_1) w_2 (x_2) = w_1 (x_1).
\eqno(5.1{\rm a})
$$
Select (as in (3.6):
$$
v_0 (x) =  v_0 (x_1) \; = \; v_1(x_1) + v_2(x_2) = 
-\frac{w_{1x_1x_1} (x_1)}{w_1 (x_1)} \; - \;  \frac{w_{2x_2x_2} (x_2)}{w_2 (x_2)}
$$
$$
= 
-\frac{w_{1x_1x_1} (x_1)}{w_1 (x_1)} +0, x \neq 0, x_{11}, 1, \; v_0 \in L^\infty (0,1).
\eqno(5.1{\rm b})$$
This  can be achieved if, for example,  we select  $w$ such that is linear near the boundary and $x_{11}$, where it vanishes.

\bigskip
Note that $ w_1 (x_1)$, satisfying to (5.1b), also  solves the following one dimensional Dirichlet spectral problem in  $ H^2 (0,1) \bigcap H_0^1 (0,1)$ ({\em with simple eigenvalues)}:
$$
\omega_{{1k}_{x_1x_1}} (x_1) + v_0 (x_1) \omega_{1k} (x_1) = \lambda_{1k} \omega_{1k} (x_1), \;\; k = 1, \ldots, \;\; \lambda_{11} > \lambda_{12}  > \ldots 
\eqno(5.2)$$
$$
\parallel  \omega_{1k}   \parallel_{L^2 (0,1)} \; = \;1, \;\; k= 1, \ldots
$$
and for some $ k = k^*$ we have
$$ 
\lambda_{1k^*} = 0, \;\;  \omega_{1k^*}  = w_1 =w, \;  \lambda_k > 0, k = 1, \ldots, k^* - 1.
$$

\bigskip

The  choice of $w_1$ in (5.1a-b) implies that 
$$
w(x_{11}, x_2) = w_1 (x_{11}) = \omega_{1k^*} (x_{11})  = 0, 
$$
and  thus, without loss of generality,  $ \omega_{1k^*} (x_1) $  is positive on the left of $x_{10}$ and negative on the right of it (it also vanishes at $x_1 = 0, 1$).

\bigskip

{\bf Remark 5.1: Instrumental properties of  solutions to   Sturm-Liouville problem.}  {\em Let us recall along these lines (\cite{Wal}, page 272)  that all eigenvalues of the Sturm-Liouville problem of our interest below, namely:
$$
z_{i x_1 x_1} + q(x_1) z_i (x_1) = \beta_i z_i (x_1)\;\; x \in (0, 1), t>0, z _i |_{x_1 = 0, 1} = 0, \; i = 1, \ldots,
$$
where $ \beta_1 < \beta_2 < \ldots$, $q \in C[0,1]$,
 can have only finitely many zeros and $z_i (x_1) $ has exactly $ i -1$ zeros. Furthermore,  between two successive zeros of  $ z_i $ 
and also between $ x_1 =0 $ and its  first zero and between its  last zero and $ x_1 = 1 $ there
is exactly one zero of $ z_{i+1}$.}

\bigskip

Due to  Remark 5.1,  {\em $ \omega_{12} (x_1) = w_1 (x_1) = w (x)$ is the 2nd eigenfunction for the spectral problem  (5.2)} (indeed,  the first eigenfunction does not have zeros in $(0, 1)$ and the next eigenfunction has one zero in $(0,1)$), and 
$$
\lambda_{1k^*} = \lambda_{12}.
\eqno(5.3) $$

\bigskip

\underline{{\em The form of eigenelements of (1.1).}}
The eigenvalues and orthonormalized  (in $L^2 (0, 1)$) eigenfunctions of the respective spectral Dirichlet problem for (1.1) with 
$ v(x_1,x_2) = v_0 (x_1),$  namely: 
$$
\Delta \omega_l (x) + v_0 (x_1) \omega_l (x) = \lambda_l \omega_l (x),
\eqno(5.4)$$
will be of the following form:
$$
\omega_l (x) = \omega_{1k} (x_1) \sqrt{2} \sin \pi m x_2, \;\;\;\; \lambda_l = \lambda_{1k} - (\pi m)^2, \; l, k, m = 1, \ldots, \;\; \lambda_1 \geq \lambda_2 \geq \ldots 
\eqno(5.5)$$

\bigskip

{\bf Step 2: Auxiliary steering to $ u_*$.} 
In terms of the momentum problem in (3.3a-b), we want to have:
$$
 \omega_{k_*} (x)  = \sqrt{2} \sin  \pi x_2 \, \omega_{12} (x_1), \;\; \lambda_{k_*} = \lambda_{12} - (\pi)^2 = - (\pi)^2.
 \eqno(5.6)$$
 
Due to (5.5) and Remark 5.2, 
\begin{itemize}
\item
$\omega_{k_*}$ is the only eigenfunction from (5.4) which has a single, more precisely, {\em vertical} line of change of sign and no horizontal lines of change of sign. 
\item
In turn,  there is only one eigenfunction in  (5.4) with a single {\em horizontal} zero-line, namely:
$$ 
\sqrt{2} \omega_{11} (x_1) \sin 2 \pi x_2,
$$
\item
while $\omega_1(x) = \omega_{11}(x_1) \sin \pi x_2$ does not have internal zero lines in $\Omega$.
\end{itemize}
 
\underline{{\em We intend, in this Step 2,}} to apply the method of Theorem 3.1 to approximately steer system (1.1), on some $ (0, T_*)$,  to a state $u |_{t=T_*} = u_* (x)$  that ``approximately'' solves  the momentum problem in (3.3a-b) for $k_*$  as in (5.6), that is, with control $v$ as in (5.1). We will to show below that it is possible. 

\bigskip
\underline{{\em Description of desirable $u^*$ in (3.5a) and $u_*$.}} If the aforementioned desirable  $u_*$ is chosen as the new initial condition for (1.1), (5.1) on $(T_*, \infty)$,  it will  ``approximately'' eliminate (on some $(T_*, T^*)$) all the terms in the respective generalized Fourier series representation (5.7) of solution to (1.1), (5.1), $u |_{t=T_*} = u_*$,   preceding (in order of decrease of $ \lambda_l$'s) the term containing the eigenfunction $ \omega_{k_*}$, namely: 
$$
u (x, t) \; = \sum_{m = 1}^{m_*} \;   \left( \int_\Omega u_* \omega_{11} (x_1) \sqrt{2} \sin \pi m x_2  dx \right) e^{(-(\pi m)^2 + \lambda_{11}) (t-T_*)} \omega_{11} (x_1) \sqrt{2} \sin \pi m x_2
 $$
 $$
  + \; 
 \left( \int_\Omega u_* \omega_{k_*} dx \right) e^{-\pi^2 (t-T_*)} \omega_{k_*} (x) \; + \; \ldots, \;\; t >T_*
 \eqno(5.7)$$
for some $ m_* \geq 1$ such  (if exists) that 
$$
-(\pi m)^2 + \lambda_{11} = \lambda_l > \lambda_{k_*} =-(\pi)^2 + \lambda_{12}= - \pi^2, l = 1, \ldots, l_* -1.
$$
The other terms, not explicitly present  in (5.7), in view of (5.5)
are associated with eigenvalues strictly smaller that  $ \; -(\pi)^2$. 

Therefore, if $u_* $ is such that it  ``approximately eliminates''  the  1st sum on the right of (5.7), that is, if, e.g., 
$$
 \int_\Omega u_* (x_1, x_2) \omega_{11} (x_1) dx_1
 $$
can be made ``as small as we wish'' in $ L^2 (0,1)$ by applying Theorem 3.1, while satisfying a  respective condition of type (3.3b), 
then we will be in a position  to apply Theorem 3.2 to steer system (1.1), (5.1)  to a term in (5.7) containing 
$ \omega_{k_*}$ at some $ t= T^*$.

Ideally, this would be the case, if, for example,  $ u_* \in L^2 (\Omega) $ is such that it has the same zero-line as $ u_0$ (this can be achieved by Theorem 3.3) and, simultaneously, in notations of (3.5b)
$$
u_* (x) = u^* (x) = u_1 (x_1) u_2(x_2), \; u_2 (x_2) \equiv 1, \;\;  \int_\Omega u_1  \omega_{11}  dx_1  \; = \; 0, \;\;\;\; \int_\Omega u_1  \omega_{12} dx _1 \; = c^0 \neq 0,
$$
or, which is the same, $ u_*$ formally solves  the $k_*$-momentum problem (3.3a-b) for our $ v_0$ in (5.1).
Unfortunately, Theorem 3.1 deals with approximate steering in $ L^2 (\Omega)$ only.  

Nonetheless, 
we can use  Theorem 3.1  to  get rid of the 1st sum in (5.7) with any {\em pre-assigned accuracy}.

\underline{{\em Selection of $u^*(x)$ as in (3.5a) and  $u_*$.}} Indeed, since $ \omega_{11} (x_1)$ does not change sign in $ (0, 1)$, and $ \omega_{12}  (x_1)$ admits only one change of sign in $(0,1)$,
we can easily select a function $ u_1  \in L^2 (0,1) $ with the same order of change of sign along the $x_1$-axis as $u_0$ in (1.1)  such that
$$
\int_\Omega  u_1 (x_1) \, \omega_{11} (x_1) dx  \; = \; 0, \;\;\;\; \int_\Omega u_1 (x_1)  \, \omega_{12} (x_1) dx \; = c^0  \; > \; 0.
\eqno(5.8)$$
In terms of (3.5a), we set
$$
u^*(x) = u_1(x_1) \times 1 = w_1(x_1) w_2 (x_2).
$$

Let us consider any  sequence $\{T_{*i}\}_{i= 1}^\infty$ (it will be further defined in more detail) such that
$$
 \lim_{i \rightarrow 0} T_{*i} = 0.
$$ 
Theorem 3.1 allows us to steer system (1.1) on each of the  intervals $(0, T_{i*})$ to  states  $  u_{i*} = u (\cdot, T_{i*}) $, with the same single line of change of sign as the original $ u_0$, such that
desirable $ u_{i*}$'s admit the following presentation:
$$
u_{i*} (x)  = u(x, T_{i*}) \; = u^* (x) + r(x, T_{i*}) \; = u_1 (x_1) + r(x, T_{i*}),
\eqno(5.8^\prime)$$ 
where 
$$
\parallel r(\cdot, T_{i*}) \parallel_{L^2 (\Omega)} \; \; \rightarrow \; 0 \;\;{\rm as} \; T_{i*} \rightarrow 0+.
 \eqno(5.9)$$

For $ t >T_{i*}$, in view of (5.7),  {\em  if we will use control } $ v(x) = v_0 (x_1)$, the respective solution to system (1.1) will have the following representation:
 $$
u (x, t) \; = \;
 \left( \int_\Omega u_1 (x_1) \omega_{k_*}  dx \right) e^{-\pi^2 (t -T_{i*})} (\int_0^1 \sqrt{2} \sin \pi x_2  dx_2)  \omega_{k_*} (x)\; 
 $$
 $$
 + \sum_{\lambda_l < -(\pi)^2}  \left( \int_\Omega u_1 (x_1)  \omega_{l} (x) dx \right) e^{\lambda_l (t- T_{i*})}   \omega_{l} (x)   \; + \; p (x, t-T_{i*}), \;\;\;\; t > T_{i*},
 \eqno(5.10{\rm a})$$

 where $ \omega_l$'s do not contain $\omega_{11}$, due to (5.8),  and $ \omega_l$'s and $ \lambda_l$'s are defined in (5.5) and 
 $$
 p (x, t- T_{i*})\; = \; \sum_{l = 1}^\infty  \left( \int_\Omega r(x, T_{i*})  \omega_{l} (x) dx \right) e^{\lambda_l (t- T_{i*})}   \omega_{l} (x),
 \eqno(5.10{\rm b}) $$
$$
 \parallel p (\cdot , t - T_{i*}) \parallel_{L^2 (\Omega)} \; \leq  \; e^{\lambda_1 (t- T_{i*})} \parallel r(\cdot, T_{i*}) \parallel_{L^2 (\Omega)} .
 \eqno(5.10{\rm c})$$

\bigskip
{\bf Step 3: Steering to $\omega_{k_*}$.}  Now we apply the argument of Theorem 3.2 in Section 4 for 
$ t  \in  (T_{i*}, T_i)$, where $\{T_i\}_{i = 1}^\infty $ is any monotone increasing sequence such that $ T_i > T_{i*}$ and 
$$
 \lim_{i \rightarrow \infty} T_i  \; = \; \infty,
 $$
 with control
 $$ 
v (x) =   \; v_0 (x_1) - \lambda_k{_*}  + a_i=   \; v_0 (x_1) + \pi^2  + a_i , i = 1, \ldots, 
$$  and with 
$$
a_i = \frac{1}{T_i - T_{i*}} \ln (\frac{1}{\int_\Omega u_1 (x_1)  ) \omega_{k_*}  dx}) \; \rightarrow  0 \;\; {\rm as} \;\; i \rightarrow \infty, 
 \eqno(5.11)$$
see (5.1a-b) and (4.4).
 Note that (5.11) implies that 
 $$
\lambda_l - \lambda_{k_*}  + a_i  < 0, \;\; \lambda_l  < \lambda_{k_*} \;\;\;\;  {\rm when}   \;\;\;\;  T_i  - T_{i*} \; \rightarrow  \; \infty  \;\; {\rm as} \; i \rightarrow \infty.
 \eqno(5.12)$$ 
 
This will result in the following formula in place of (4.2) and (5.10a-c):
$$
u (x, T_i) \; = \; \omega_{k_*} (x) \; + \; \sum_{\lambda_l < -\pi^2}  \left( \int_\Omega u_1 (x_1)  \omega_l dx \right) e^{(\lambda_l -\lambda_{k_*}  + a_i) (T_i- T_{i*})} \omega_l(x)
$$
$$
+ \; \sum_{l = 1}^\infty  \left( \int_\Omega r(x, T_{i*})  \omega_l dx \right) e^{(\lambda_l -\lambda_{k_*}  + a_i) (T_i - T^*_{i*})} \omega_l(x), \;\; i = 1, \ldots,
\eqno(5.13)$$
 where, again  $ \omega_l$'s do not contain $\omega_{11}$, due to (5.8).
 
In view of (5.11)-(5.12), the 1st  series on the right in (5.13) tends to zero in $ L^2 (\Omega)$   as $ T_i - T_{i*}  \rightarrow \infty$.

The same, by the same reasoning, will happen to the tail of the 2nd  series beginning for $ \lambda_l < -\pi^2$, see 
(5.9). In turn, for finitely many other terms for the 2nd series in (5.13) (namely, $k_* - 1$ terms) we have the following estimate:
$$
||  \sum_{\lambda_l >-\pi^2}  \left( \int_\Omega r(x, T_{i*})  \omega_l dx \right) e^{(\lambda_l -\lambda_{k_*}  + a_i) (T_i - T^*_{i*})} \omega_l ||_{L^2 (\Omega)}
$$
$$
\leq (k_*-1)  e^{(\lambda_1 -\lambda_{k_*}  + a_i) T_i}  \parallel r(\cdot, T_{*i}) \parallel_{L^2 (\Omega)} \; \rightarrow 0 \;\;{\rm as} \;\; i \rightarrow \infty, 
$$
if (for a selected sequence of $T_i$'s) $T_{i*}$'s are, additionally to their convergence to zero, selected   to ensure the rate of convergence to zero in (5.9) to be higher than that of convergence of $ e^{(\lambda_l -\lambda_{k_*}  + a_i)T_i}$ to $\infty$  as $ i \rightarrow \infty$. 

Combining the above yields that  $ u (\cdot, T_i ) \rightarrow  \; \omega_{k_*} (\cdot)$ in $ L^2 (\Omega)$ as $ i \rightarrow \infty$. 

\bigskip
{\bf Step 4: ``Magnitude adjustment'' steering.} 
This can be achieved by  applying Theorem 3.1 to system (1.1)  on some  intervals $(T_i, T^{(i)}) $ to steer it to a sequence of states converging to $ u_1$.

This ends the proof of Theorems 3.3-3.4 in the  case of squared domain and a single line of change of sign. 
\bigskip


{\bf 6. Proof of Theorems 3.3 and 3.4  in the general case: ``Separation of variables''.}  To this end, we only need to modify the arguments in Steps 1-2 from Section 5, while the Steps 3 and 4 would be nearly identical. Namely,  we need to achieve an auxiliary steering, on some $(0, T_*)$, from the original initial state  to an intermediate state $u_* (x)$  with the same zero-hyperplanes as $ u_0$ such that  the the term containing the desirable  $ \omega_{k_*} $   (that is, with the same zero-hyperplanes as the target state $ u_1(x)$) will be the 1st ``substantially nonzero'' term in the respective expansion of  solution  like in (5.7) after $T_*$ when  a suitable control applied . 

\bigskip

Without loss of generality and for simplicity of notations, we can assume that $\Omega$ is a unit $n$-dimensional cube.

\bigskip

{\bf Separation of variables.}  Let us assume that our target state $ u_1 (x)$ has $k_* -1$ hyperplanes of chage of sign $P_{ij}, i = 1, \dots, n, j = 0, \ldots, k_i -1$ as given in (3.4b).

Select  functions  $u_i (x_i)$'s, $u^*$,  $ w_{i} (x_i)$'s, $w$, and control $v_0$ as in (3.4a) -(3.6).
In turn,  in this general case, instead of (5.5)-(5.6), we have (see (3.7a) for notations):
$$
\omega_{k_*} (x) =    \Pi_{i=1}^n   \omega _{i k_i} (x_i).
\eqno(6.1)$$  

Respectively, to make the term with $ \omega_{k_*} (x)$ be the first (in the order of decrease of $ \lambda_i$'s) essentially non-zero terms, making use of Theorem 3.1, we need to steer (1.1) on some $ (0, T_*)$ ($T_* \rightarrow 0+ $  as in Step 2 in Section 5)  to
$$
u_* (x) = \Pi_{i=1, \ldots, n, j_i = 1, \ldots}    u_{i} (x_i),
\eqno(6.2)$$
which has the same zero-hyperplanes as the initial state $ u_0$, such that 
 (3.7b)-(3.8b) holds, as it is assumed in Theorem 3.4. This can be achieved by dealing separately with each $i$-th one dimensional problem, $i = 1, \ldots, n$ (similar to how it was done  in Section 5 for $i = 1$) and with the same $T_* (\rightarrow 0+ )$ for all $i$'s.  Note that the latter conditions are straightforward for Theorem 3.3 as we discussed it in Section 5.

This ends the proof of Theorems 3.3 and 3.4.



\bigskip
{\bf 7. Proof of Theorems 3.5 and 3.6.}  Let us show  how one can find $u_{i} (x_i)$'s,  satisfying  (3.7b)-(3.8b).

\bigskip

{\bf 7.1. Conversion of (3.7b)-(3.8b) to a problem associated with a system of  linear algebraic equations.} 
For simplicity of notation assume again that $ \Omega = (0, 1) \times \ldots \times (0, 1)$.

Pick any $ i \in \{1, \dots, n\}$. Without loss of generality we can further  assume that $ i = 1$ and   $k_1-1 \geq 1$, that is, there is at least one zero-hyperplane perpendicular to the $x_1$-axis. Select one more point $s$ in $(0,1)$  different from $x^0_{1j}$'s in (3.4a) (we will refine this selection later).

Next, select  $ u_{1}  (x_{1}) $  in (3.4a) to be a piecewise constant function defined by  $2k_1 +1$ real values  $S, l_j, m_j, j = 1, \ldots, k_1 $ as follows:
$$
u_{1} (x_1)  = 
\left\{ \begin{array}{ll}
l_j, \;\;\;\; &  x_1 \in (x^0_{1j} -h, x^0_{1j}),   \\ 
m_j, \;\;\;\; &  x_1 \in (x^0_{1j}, x^0_{1j}+ h),   \\ 
S, \;\;\;\; &  x_1 \in (s, s+ h),   \\ 0, \;\;\;\; & {\rm elsewhere \;  in}\;  (0, 1).
\end{array}
\right.
\eqno(7.1)$$
This $u_1 (x_1)$ can be reached for (1.1) as close as we wish in $ L^2 (0,1)$ due to Theorem 3.1.

{\bf Remark 7.1: Values of $ (l_j + m_j) $'s and of $S$.} {\em Note that $l_j m_j \leq 0, j = \ldots, k_1$ as $x^0_{1j}$'s are points of change of sign for $ u_1 (x_1)$. The signs of $S,  l_j$'s and $ m_j$'s are defined by the sequence of change of sign of  $ u_1 (x_1)$, or, which is the same, of $u_0 (x)$ along the $x_1$-axis.  Hence, the range of available $ (l_j + m_j) $'s in (7.1) is $R$, while $ S$ cannot change its sign.}

\bigskip
\underline{Our goal below} is to investigate (7.1))  when $ h $ tends to zero. Therefore, without loss of generality, we can assume that intervals  $ (x^0_{1j}-h_j, x^0_{1j} + h_j) $  do not overlap.

 \bigskip

The use of (7.1) will generate the following linear algebraic system of $k_1 -1$ equations in  $ l_j, m_j  $ out of conditions (3.7b)-(3.7c) in Theorem 3.4:

$$
\sum_{j=1}^{k_1 -1} A_{kj} l_j  + \sum_{j=1}^{k_1  -1} B_{kj} m_j  + C_k S=0, \;\; k = 1, \ldots, k_1 -1,
\eqno(7.2)$$
where
$$
A_{kj} = \int_{x^0_{1j} -h}^{x^0_{1j}} \omega_{1k} (x_1) dx_1, 
$$
$$
B_{kj} = \int_{x^0_{1j}}^{x^0_{1j}+h}   \omega_{1k} (x_1) dx_1, 
$$
$$
C_{k} = \int_s^{s+h}   \omega_{1k} (x_1) dx_1, \;  k= 1, \ldots, k_1 -1.
$$

Since 
$$
 \omega_{1k} (x_1) = \omega_{1k} (x^0_{1j}) + \omega^\prime _{1k} (x^0_{1j}) ( x_1 - x^0_{1j}) +O ((x_1 - x^0_{1j})^2)
$$
 in $ (x^0_{1j}-h, x^0_{1j}+ h) $ as $h \rightarrow 0$, and 
 $$
 \omega_{1k} (x_1) = \omega_{1k} (s) + \omega^\prime _{1k} (s) ( x_1 - x^0_{1j}) + O ((x_1 -s_j)^2)
$$
in $ (s, s + h)$, $ k = 1, \ldots, k_{1} -1$,
we can re-write (7.2) as follows:
$$
h \sum_{j=1}^{k_1 -1}  (l_j + m_j) \omega_{1k} (x^0_{1j}) +  h S  \omega_{1k} (s) 
= + O(h^2) [ \max  |(l_j |  + | m_j| ) + | S |], \;\; j = 1, \ldots, k_1-1.
\eqno(7.3)$$ 

\bigskip
{\bf An auxiliary linear algebraic system on a \underline{cone}.}  Let us consider the following ``limit'' system (7.4) of $(k_1-1) $ linear algebraic equations in $k_1$ real-valued variables $(V_1, \ldots, V_{k_1-1}, P)$:
\begin{itemize}
\item
 $P$ of certain sign only, defined by the location of $s$ between $x^0_{1j}$'s and the sequence of change of sign of $u_0$ along the $x_1$-axis,
 \item
 and $V_j  \in R, j= 1, \ldots, k_{1} -1$ (see Remark 7.1), 
\end{itemize}
generated by (7.3) as $h \rightarrow 0$:
 $$
\sum_{j=1}^{k_1 -1}  V_j \omega_{1k} (x^0_{1j}) +  P \omega_{1k} (s) = 0, \;\; k = 1, \ldots, k_1-1,
\eqno(7.4)$$
Since the number of equations exceeds the number of variables, without loss of generality, we can say that (7.4) admits  non-trivial solutions under the above restrictions of variables forming a cone in $R^{k_1}$. Furthermore, under Assumption 3.1 we can consider $P$ to be a free parameter, while under Assumption 3.2 we can set $ P=0$.

Indeed, we can solve this system in $ R^{k_1}$, and then if it gives solution $(V_1, \ldots, V_{k_1-1}, P)$ with $P$ of a ``wrong'' sign, we can replace (7.4) with:
 $$
\sum_{j=1}^{k_1 -1}  V_j^* \omega_{1k} (x^0_{1j}) +  P^*\omega_{1k} (s) = 0, \;\; k = 1, \ldots, k_1-1,
\eqno(7.4)^\prime)$$
and further deal with solution $(V_1^*=-V_1, \ldots, V_{k_1-1} =-V_{k_1-1}, P^* = -P)$, where  $ P^* = -P$ will be of a ``correct'' sign.

\bigskip
{\bf Construction of $u_1^h(x_1)$ in (3.7c).}  Consider any solution  $(V_1, \ldots, V_{k_1-1}, P)$ to (7.4) with ``correct sign'' of $P$ and {\bf fix it}. Then set:
$$
S =\frac{P}{h},  V_{j}  = \frac{l_{j}}{h}, m_j =0 \;\; {\rm or} \;\;  V_{j}  = \frac{m_{j}}{h}, l_j = 0\;\;  j= 1, \ldots, k_{1} -1,
$$
depending on the sign of $V_j$'s. 

Construct next $ u_1 (x_1)$ as in (7.1) with the just described choice of $ S$ and $l_{j} $'s and $ m_{j}$'s.
This will give us, due to (7.3), (7.4) the relation in  (3.7c) for $u_1^h$  in the following form:
$$
 \int_0^1 u_1^h  (x_1)  \omega_{1j} (x_i) dx_1  \; = \; O(h),  \; j = 1, \ldots, k_1 -1 \;\; {\rm as} \; \; h \rightarrow 0.
 \eqno(7.5) $$
 
 \bigskip
 {\bf Remark 7.2} {\em Note that relation (7.5) will hold uniformly of any aforementioned  $(V_1, \ldots, V_{k_1-1}, P)$ lying in a fixed bounded set, e.g., in a fixed ball intrersected with a cone describing the restriction on sign of $S$.} 
 
 \bigskip
{\bf Discussion of  condition in (3.8b).}  To ensure a suitable  ``contribution of $ u_1^h (x_1)$'s in (3.8b), it suffices to show that we can ensure that we can find a sequence $\{h_l\}_{l=1}^\infty$ such that for $u_1^{h_l}   (x_1) $ satisfying (7.5) we can have:
$$
|  \int_0^1 u^{h_l}_{1} (x_i)  \omega_{1k_1} (x_1) dx_1  | \;  = 1,  \;\;{\rm as} \; h_l \rightarrow 0, \; l = 1, \ldots
\eqno(7.6)$$

 \bigskip
In view of Remark 7.1, to achieve  (7.6), it is sufficient to show that vector
$$
(\omega_{1k_{1}} (x^0_{11}), \ldots, \omega_{1k_{1}} (x^0_{1k_1-1}),  \omega_{1k_{1}} (s))
\eqno(7.7)$$ 
does not belong to the span of vectors
$$
(\omega_{1k} (x^0_{11}), \ldots, \omega_{1k} (x^0_{1k_1-1}),  \omega_{1k} (s)), \;\; k \leq k_1-1,
\eqno(7.8{\rm a})$$ 
that is, 
$$
(\omega_{1k_{1}} (x^0_{11}), \ldots, \omega_{1k_{1}} (x^0_{1k_1-1}),  \omega_{1k_{1}} (s)) \; \not\in \;\;{\rm span} \; \{(\omega_{1k} (x^0_{11}), \ldots, \omega_{1k} (x^0_{1k_1-1}),  \omega_{1k} (s)) \; | \;\; k \leq k_1-1\}.
\eqno(7.8{\rm b})$$ 
If this is true for the ``cut-off'' vectors composed of the vectors in (7.7) and (7.8), by removing the the last coordinates, as in Assumption 3.2, then we are done.

If Assumption 3.2 does not hold, then, under Assumption 3.1, the vector in (7.7) would be a unique linear  combination of the vectors in (7.8) because the  vectors 
$$
 (\omega_{1k} (x^0_{11}), \ldots, \omega_{1k} (x^0_{1k_1-1})), \; k = 1, \ldots,  k_1-1,
 \eqno(7.9)$$
composed of  the first $ k_{1}-1$ coordinates in (7.8), are linear independent.

Recall now that we are free to select a point $s$ to avoid (7.8b). Indeed, if this is not possible for any $ s$ in $(0,1)$, then the last coordinate in (7.7) or (7.8), namely, $\omega_{1k_{1}} (x_1)$ will be a unique linear combination of the functions 
$ \omega_{1k} (x_1),   k =1, \ldots, k_{1}-1$ when $x_1$ ranges over  all $(0,1)$) (see the previous paragraph). But $ \omega_{1k} (x_1),   k =1, \ldots, k_1$ are orthonormal to each other in $ L^2 (0,1)$ and, hence, are linear independent. Moreover, the set of such ``bad'' points $s$ for which (7.8b) dos not hold,  can at most be countable  as we discussed it after Assumptions 3.1 and 3.2 in Section 3.

This completes the proof of Theorems 3.5 and 3.6.

\bigskip

\bigskip
{\bf Appendix:  Proof of Theorem 3.1.}  We intend to adapt the scheme, previously used for system (1.1) in \cite{Kh3} (Chapter 3)  in several spatial dimension, assuming that its solutions do not change  sign,  or in \cite{CanKh2} in the case in one spatial  dimension, assuming that system's solutions  can change sign finitely many times.  Our goal here is to apply a suitable static control $ v = v(x)$ which will steer system (1.1) from a given $u_0$ to any desirable  state $ u (\cdot, T)$ (satisfying conditions stated in Theorem 3.1)  as close in $ L^2 (\Omega)$ to $ u_1$  as we wish. We   split the proof into  two cases.

\bigskip
{\bf Step 1. A special case.}
In this step we assume the following condition:

\bigskip

{\bf Assumption A.1.} {\em Suppose that $ \mid u_1 (x)  \mid \, < \, \mid u_0 (x) \mid $ in the  sets $ S_i^+$'s and $ S_i^-$'s  in Definiton  3.1, and 
the function
$$
v_0 (x) \; = \;  \left\{ \begin{array}{ll}
\ln \left( \frac{u_{1} (x)}{u_0 (x) } \right), \;\;\;\; &  {\rm where} \;\; u_0 (x) \neq 0 \; {\rm in} \; \Omega,   \\ 
0, \;\;\;\; &  {\rm elsewhere} \; {\rm in} \;  \Omega, \; {\rm that \; is, \; on \; some \; sets \;of \; measure \;  zero}\\
\end{array}
\right.   
\eqno(A.1)$$ 
lies in $ L^\infty (\Omega)$. 
}

\bigskip
{\bf Step 1.1: Selection of bilinear control.}
Note that the function $\, v_0 (x) $ in (A.1) satisfies 
$
v_0 (x) \; \leq \; 0 $ in $ \Omega$.
Let 
$$
v (x,t) \; = \; \frac{1}{T} v_0 (x).
\eqno(A.2)$$
Then, the corresponding solution to  (1.1), treated as an ordinary differential equation  in time in a respective Banach space,  admits the following representation a.e. in $\Omega$:
$$
u(x,t) \; = \; e^{v_0 (x) \frac{t}{T}} u_0 (x)  \; + \; \mathop{\int}_0^t e^{v_0 (x) \frac{(t - \tau)}{T}} \Delta u (x, \tau) d \tau\,.
$$
At time $ t = T$ we have a.e. in $\Omega$:
$$
u(x,T) \; = \; u_1 (x)  \; + \; \mathop{\int}_0^T e^{v_0 (x) \frac{(T - t)}{T}} \Delta u (x, t) dt.
\eqno(A.3)$$

\bigskip
{\bf Step 1.2: Evaluation of $ \parallel  \Delta u \parallel_{L^2 (Q_T)}$.} Let us show that the 2-nd term in the right-hand side of (A.3) tends to zero in $ L^2 (\Omega) $ as $ T \rightarrow 0+$, which would mean that $ u(\cdot, T) \rightarrow u_1 $ in $ L^2 (\Omega)$ at the same time.
Note first that, since $ v_0 (x) $ is nonpositive,
$$
\mathop{\int}_\Omega \left( \mathop{\int}_0^T e^{v_{0} (x)
\frac{(T - \tau)}{T}} \Delta u (x, \tau) d \tau \right)^2 dx \; \leq \;
T \parallel  \Delta u \parallel^2_{L^2 (Q_T)}.
\eqno(A.4)$$

Without loss of generality, we can further assume that $ v_0 \in C^2(\bar{\Omega})$. 

\bigskip
{\bf Remark A.1.}
 {\em Indeed, if $ v_0 \not\in  C^2(\bar{\Omega})$, then we could consider instead a sequence of uniformly bounded controls $ \{ v_{0l}\}_{l = 1}^{\infty},  \, v_{0l} \in C^2(\bar{\Omega})$, approximating $ v_0 $ in $ L^2 (\Omega) $, making use of the following limit relation:
$$
e^{v_{0l} (x) t/T} u_0 (x)  \mid_{t = T} \; \rightarrow  \;
e^{v_0 (x) t/T} u_0 (x) \mid_{t = T} \; = \; u_1 (x) \;\;{\rm in} \;\; L^2 (\Omega) \;\; {\rm as} \; l \rightarrow \infty.
$$
}

\bigskip 
Multiplying (1.1) by $u_{xx}$ with $ v = \frac{1}{T} v_{0} \leq 0 $ and  integrating by parts over $ Q_T$, we have:
$$
\parallel  \Delta u  \parallel^2_{L^2 (\Omega \times (0,T))} \; = \;
\mathop{\int}_0^T \mathop{\int}_\Omega u_{t} \Delta u \, dx dt \; - \;
\frac{1}{T} \mathop{\int}_0^T  \mathop{\int}_\Omega v_{0} u \Delta u \, dx dt
$$
$$
= \; - \frac{1}{2} \mathop{\int}_0^T \mathop{\int}_\Omega  \parallel \nabla u \parallel^2 _t dx dt \; + \;
\frac{1}{2T} \mathop{\int}_0^T  \mathop{\int}_\Omega \nabla v_{0} \cdot  \nabla \parallel u \parallel^2 dx dt \; + \;
\frac{1}{T} \mathop{\int}_0^T \mathop{\int}_\Omega v_{0}  \parallel \nabla u \parallel^2 dx dt
$$
$$
= \; -  \; \frac{1}{2} \mathop{\int}_0^T \mathop{\int}_\Omega \parallel \nabla u \parallel^2 _t dx dt 
 \; - \;
\frac{1}{2T} \mathop{\int}_0^T  \mathop{\int}_\Omega \Delta v_0  \, u^2 dx dt\; +  \frac{1}{T} \mathop{\int}_0^T \mathop{\int}_\Omega v_{0} \parallel \nabla u \parallel^2  dx dt.
$$
Thus, we obtain:
$$
\parallel \Delta  u \parallel^2_{L^2 (\Omega \times (0,T))} \; +  \; \frac{1}{2} \mathop{\int}_\Omega  \parallel \nabla u (x,t) \parallel^2  dx  \; - \; \frac{1}{T} \mathop{\int}_0^T \mathop{\int}_\Omega v_{0} \parallel \nabla u \parallel^2  dx dt
$$
$$
= \; \frac{1}{2} \mathop{\int}_\Omega  \parallel \nabla u (x,0) \parallel^2  dx \; - \;
\frac{1}{2T} \mathop{\int}_0^T \mathop{\int}_\Omega \Delta  v_{0} \, u^2 dx dt.
\eqno(A.5)$$
In particular, recalling that $ v_0 (x) \leq 0 $,
$$
\parallel \Delta  u \parallel^2_{L^2 (\Omega \times (0,T))}  \; \leq  \;
\frac{1}{2} \mathop{\int}_\Omega  \parallel \nabla u_0 \parallel^2   dx \; + \;
\frac{1}{2T} \max_{x \in \bar{\Omega}} \mid \Delta v_{0} \mid
\mathop{\int}_0^T \mathop{\int}_\Omega   u^2 dx dt.
$$

Now, since $ v_0 (x) \leq 0$, multiplication of (1.1) by $ u$ and integration by parts over $ Q_T$ yield: 
$$
\mathop{\int}_0^T \mathop{\int}_\Omega u^2 (x,t) dx dt\; \leq \; T   \mathop{\int}_\Omega u_0^2 (x) dx.
$$
Hence, combining this with (A.5), we derive that
$$
\parallel  \Delta u \parallel^2_{L^2 (\Omega \times (0,T))} \; \leq  \;
\frac{1}{2} \mathop{\int}_\Omega  \parallel \nabla u_0 \parallel^2   dx \; + \;
\frac{1}{2} \max_{x \in \bar{\Omega}} \mid \Delta v_{0} \mid
\mathop{\int}_0^T \mathop{\int}_\Omega   u_0^2 dx.
\eqno(A.6)
$$
In turn,  making use of (A.4) and (A.6), we further obtain that
$$
\mathop{\int}_\Omega \left( \mathop{\int}_0^T e^{v_{0} (x)
\frac{(T - \tau)}{T}} \Delta u (x, \tau) d \tau \right)^2 dx \;
\leq \;  \frac{T}{2} \mathop{\int}_\Omega  \parallel \nabla u_0 \parallel^2   dx \; + \;
\frac{T}{2} \max_{x \in \bar{\Omega}} \mid \Delta v_{0} \mid
\mathop{\int}_0^T \mathop{\int}_\Omega   u_0^2 d.
\eqno(A.7)
$$
Note now that the right-hand side in (A.7) tends to zero as $ \; T \rightarrow 0 $, which combined with (A.3) and (A.4) yields the desirable approximate controllability  result under Assumption A.1.

\bigskip
{\bf Step 2:    Assumption A.1 does not hold}. 
In this case,  we will apply, first, an auxiliary constant  control  $ v (x,t) = m > 0$ on some interval $ (0, t_*)$, which generates the following solution to (1.1):
$$
u (x,t_*) = \;  e^{mt_*} \sum_{k = 1}^\infty  2 e^{-(\pi k)^2 t_*} \left(\mathop{\int}_0^1 u_0 (r) \omega_k ¨  \, dr \right) \omega_k (x)
$$
$$
= \; e^{mt_*} \sum_{k = 1}^\infty  2 (e^{-(\pi k)^2 t_*} -1) \left(\mathop{\int}_0^1 u_0 (r) \omega_k (r) \, dr \right) \omega_k (x)
\; + \; e^{mt_*} u_0 (x).
$$

Consider any $ L > 1$. Then, by selecting  $ m = (\ln L)/t_*, t_* > 0$,  we have that
$$
e^{mt_*} = L,
$$
and $ u (\cdot, t_*) \rightarrow L u_0 $ in $ L^2 (Q_T) $ as $ t_* \rightarrow 0+$.

Hence, for any positive integer $ i $ we can find a suitably large parameter $ L_i $ and respective moment $ t_i > 0$ such that the inequality in the first line of  Assumption A.1 will hold for $ u (x, t_i)$, regarded as a new initial condition  for a future action as in Step 1, everywhere except, possibly, some  some (measurable) set $ A_i \subset \Omega$ whose measure tends to zero as $ n$ increases.

Now, in place of the control in (A.2), we can select its modified version $ v_{0n} \in L^\infty (0,1) $ as follows:
$$
v_{0n} (x) \; = \; 
\left\{ \begin{array}{ll}
\ln \left( \frac{u_{1}(x)}{u (x, t_n)} \right), \;\;\;\; &  x \in \Omega \backslash A_i,   \\ 
0, \;\;\;\; & x \in A_i,
\end{array}
\right.    \;\;\;\; i = 1, 2, \ldots\eqno(A.8)$$

The argument of Step 1 yields (see (A.3) and (A.7))  that we can steer   the solution of (1.1) from $ u (\cdot, t_i) $ to an auxiliary target state 
$$
u_{1i} (x) \; = \; 
\left\{ \begin{array}{ll}
u_1, \;\;\;\; &  x \in \Omega \backslash A_i,   \\ 
0, \;\;\;\; & x \in A_i,
\end{array}
\right.    \;\;\;\; i = 1, 2, \ldots
$$
as close in $ L^2 ( \Omega)$   as we wish at some moment $ t_{i*} >  t_i$ (see the formula between (A.2) and (A.3), and (A.6)).
Since $ u_1 \in L^2 (Q_T)$ (and is fixed), selecting appropriately large value for $ i$ will provide that same approximate controllability result as in Step 1.

This ends the proof of Theorem 3.1.

\bigskip

\newpage


\end{document}